\documentclass[11pt]{article} 
\usepackage[latin1]{inputenc}
\usepackage[T1]{fontenc}
\usepackage[francais]{babel}
\usepackage{amsmath}
\usepackage{amsfonts} 
\usepackage{amssymb} 
\usepackage{amsthm}
\usepackage{graphicx}
\usepackage{stmaryrd} 
\usepackage[np]{numprint}
\usepackage{url}
\usepackage{vmargin}
\usepackage{mathrsfs}
\theoremstyle{plain}
\newtheorem{theoreme}{Théorème}
\newtheorem{lemme}[theoreme]{Lemme}
\newtheorem{propriete}[theoreme]{Propriété}
\newtheorem{corollaire}[theoreme]{Corollaire}
\theoremstyle{definition}
\newtheorem{definition}[theoreme]{D\'efinition}

\newtheorem{ex}[theoreme]{Exemple}

\theoremstyle{remark}
\newtheorem{rem}[theoreme]{Remarque}

\newcommand{\intervalleff }[2]{\left[{#1}\mathpunct{};{#2}\right]}

\newcommand{\intervalleoo }[2]{\left]{#1}\mathpunct{};{#2}\right[}
\newcommand{\ensemblenombre }[1]{\mathbb{#1}}
\newcommand{\N}{\ensemblenombre{N}}
\newcommand{\Z}{\ensemblenombre{Z}}
\newcommand{\Q}{\ensemblenombre{Q}}

\newcommand{\entiere}[1]{\left\lfloor #1 \right\rfloor}
\newcommand{\abs}[1]{\left\lvert#1\right\rvert}
\newcommand{\enstq}[2]{\left\{#1\,\middle|\,#2\right\}}
\setmarginsrb{2cm}{2cm}{2cm}{2.5cm}{0cm}{0cm}{0cm}{1cm}
\date{} 
\begin{document}

\title{\textbf{L'équation diophantienne \boldmath$ax^2-by^2=1$ }}

\author{Lionel Ponton \\ \small \texttt{lionel.ponton@gmail.com}}

\maketitle

\begin{abstract} On propose une méthode de résolution effective de l'équation diophantienne $(E_2): ax^2-by^2=1$ où $a$ et $b$ sont des entiers naturels non nuls et premiers entre eux. Cette méthode s'appuie sur le développement en fraction continuée de certains nombres irrationnels quadratiques que l'on décrit complètement. On commence par utiliser ces développements pour résoudre certaines équations de Pell-Fermat généralisées avant d'appliquer à l'équation $(E_2)$ les résultats obtenus. \end{abstract} 

\bigskip

\bigskip

\section{Introduction}
Un exemple classique d'application du théorème de Gauss en arithmétique est la résolution explicite de l'équation diophantienne binôme homogène du premier degré 
\[(E_1): ax-by=1\]
d'inconnue $(x,y)\in\Z^2$ où $a$ et $b$ sont des entiers relatifs non nuls et premiers entre eux. On montre sans difficulté que l'ensemble des solutions est $\enstq{(x_0+kb, y_0+ka)}{k\in\Z}$ où $(x_0,  y_0)$ est une solution particulière de l'équation $(E_1)$, solution que l'on peut toujours déterminer, par exemple, à l'aide de l'algorithme d'Euclide. Notons qu'une autre méthode de résolution (mais qui, dans le fond, utilise essentiellement les mêmes ingrédients) est d'utiliser le développement en fraction continuée du rationnel $\frac{a}{b}$ (voir \cite[p. 37]{Du07}).

Sous l'hypothèse $a$ et $b$ premiers entre eux, l'équation $(E_1)$ admet donc toujours une infinité de solutions. Il est naturel de se demander si cela reste vrai dans le cas de l'équation diophantienne binôme homogène de degré $n$ 
\[(E_n): ax^n-by^n=1\]
où $n\in\N^*$ et $a$ et $b$ sont des entiers relatifs non nuls et premiers entre eux. Pour $n\geq 3$, la réponse, négative, est une conséquence d'un théorème dû à A. Thue (voir le paragraphe \ref{secThue}). Un inconvénient majeur de ce résultat est de ne pas être effectif. Autrement dit, le théorème de Thue assure que l'équation $(E_n)$ n'a qu'un nombre fini de solutions pour $n\geq 3$ mais sa démonstration ne donne pas le nombre de solutions ni même un algorithme qui permettrait de borner ce nombre pour des valeurs de $a$ et $b$ données. L'équation $(E_n)$ a fait l'objet de nombreuses recherches au cours du XX\up{e} siècle notamment de la part d'A. Thue, C. Siegel, B. Delaunay, T. Nagell, V. Tartakovskii, W. Ljunggren, J.-H. Everste, M. Mignotte et B. W. de Weger qui ont successivement traité différentes valeurs de $n$ et obtenu des bornes générales du nombre de solutions de plus en plus précises. En 2001, en traitant les cas restants à l'aide de techniques arithmétiques très avancées, M. Bennett a finalement démontré que l'équation $(E_n)$ admet au plus une solution $(x,y)\in\N^2$ telle que $xy\neq 0$. A ce sujet et pour des références plus précises, on pourra consulter \cite{Be01}.

\bigskip

Entre $(E_1)$ qui a toujours une infinité de solutions et $(E_n)$ qui a toujours un nombre fini (éventuellement nul) de solutions pour $n\geq 3$, l'équation binôme homogène de degré 2
\[(E_2): ax^2-by^2=1\]
occupe une place à part. Si le cas particulier $a=1$, connu sous le nom d'équation de Pell-Fermat, est abondamment traité dans la littérature (voir le paragraphe \ref{secPell}), le cas général est curieusement très rarement abordé. Cet article se propose donc de présenter une méthode effective de résolution de $(E_2)$ i.e. de présenter un algorithme permettant de déterminer les valeurs de $a$ et $b$ pour lesquelles $(E_2)$ admet des solutions et, dans ces cas, de donner une description de l'ensemble des solutions.

Notre démarche reprend en la prolongeant celle exposée par D. Duverney dans \cite[pp. 40-43]{Du07}. Nous obtenons des résultats essentiellement semblables à ceux de R. A. Mollin (\cite{Mo01} et \cite{Mol04}) mais en restant à un niveau toujours élémentaire (i.e. sans faire référence aux unités des corps quadratiques). Dans le paragraphe \ref{secirrquad}, nous rappelons les résultats essentiels dont nous aurons besoin sur les développements en fraction continuée des réels avant d'étudier plus particulièrement les développements des irrationnels quadratiques et plus précisément les nombres de la forme $\entiere{\sqrt{d}}+\sqrt{d}$ où $d$ est un entier naturel qui n'est pas un carré parfait. Les résultats de ce paragraphe sont ensuite utilisés pour démontrer les deux théorèmes principaux dont le théorème \ref{theo04110947} qui donne une méthode effective de résolution de $(E_2)$.  

Nous ne prétendons pas que ce qui suit soit fondamentalement nouveau. Cependant, nous n'avons trouvé aucune référence dans la littérature existante sur le sujet qui donne une description effective de l'ensemble des solutions de $(E_2)$ telle que nous la proposons.

\section{Développement en fraction continuée des irrationnels quadratiques}
\label{secirrquad}
\subsection{Quelques rappels sur les développements en fraction continuée}

L'algorithme des fractions continuées (régulières) -- appelée plus fréquemment fractions \emph{continues} bien que cette dénomination soit malheureuse -- est un outil central en approximation diophantienne c'est-à-dire dans la recherche d'approximation de nombres réels par des nombres rationnels. En effet, il fournit de très bonnes approximations (et même, en un certain sens, les meilleures approximations possibles) et cette propriété est un argument essentiel dans divers raisonnements (voir par exemple les articles de X. Caruso dans le numéro 123-1 de la RMS ou celui de T. Lafforgue et B. Randé dans le numéro 125-1).

Différentes approches de la notion de fraction continuée sont possibles, notamment en utilisant des matrices (voir \cite[chap. II]{Fai91}) ou en adoptant un point de vue géométrique (voir à ce sujet le bel article de R. Bondil dans le numéro 124-1 de la RMS). Nous nous contentons ici de présenter l'algorithme et d'énoncer les propriétés dont nous aurons besoin par la suite. 

\bigskip

\'Etant donné un réel quelconque $a_0$, on pose $\left[a_0\right]:=a_0$ et, si $a_1$, $a_2$, ..., $a_n$ ($n\geq 1$) sont des nombres réels strictement positifs , on pose
\[ \left[a_0, a_1, ..., a_n\right] := a_0+\dfrac{1}{a_1+\dfrac{1}{a_2+\dfrac{1}{\ddots+\dfrac{1}{a_n}}}}.\]

Il est immédiat que si les nombres $a_0$, $a_1$, ..., $a_n$ sont tous entiers alors $\left[a_0, a_1, ..., a_n\right]$ est rationnel. Le lemme suivant, qui se démontre aisément par récurrence (\cite[p. 166  et suiv.]{HW60}), permet de déterminer précisément le numérateur et le dénominateur de l'écriture sous forme irréductible de ce nombre rationnel.

\begin{lemme} --- Soit $(a_n)_{n\in\N}$ une suite d'entiers, strictement positifs à partir du rang $1$, et $(p_n)_{n\geq -2}$ et $(q_n)_{n\geq -2}$ les suites d'entiers définies par:
\begin{equation}
\begin{cases} p_{-2}=0,~p_{-1}=1 \\ \forall n\in\N,~p_{n}=a_{n}p_{n-1}+p_{n-2} \end{cases} \qquad \text{et} \qquad \begin{cases} q_{-2}=1,~q_{-1}=0 \\ \forall n\in\N,~q_{n}=a_{n}q_{n-1}+q_{n-2} \end{cases}.
\label{eq29101445}
\end{equation}
\begin{enumerate}
\item La suite $(q_n)$ est croissante à partir du rang $n=-1$ et, si $a_0\geq0$, la suite $(p_n)$ est croissante à partir du rang $n=-2$.
\item Pour tout $n\in\N$, $p_{n-1}q_{n-2}-q_{n-1}p_{n-2}=(-1)^{n}$ et $q_{n}p_{n-2}-p_{n}q_{n-2}=(-1)^{n-1}a_{n}$.
\item Pour tout $n\in\N$ et pour tout réel $y>0$, $\left[a_0, a_1, ..., a_{n}, y\right]=\dfrac{p_{n}y+p_{n-1}}{q_{n}y+q_{n-1}}$.
\item Pour tout $n\in\N$, $\dfrac{p_n}{q_n}$ est l'écriture sous forme de fraction irréductible de $\left[a_0, a_1, ..., a_n\right]$.
\end{enumerate}
\label{lem29101626}
\end{lemme}

\bigskip

\noindent \textbf{Notation} --- Dans toute la suite, si $x$ est un réel, on note $\entiere{x}$ la partie entière de $x$ i.e. l'unique entier $k$ tel que $k \leqslant x < k+1$.

\begin{definition} --- Soit $x$ un nombre réel quelconque. On associe à $x$ deux suites (éventuellement finies): une suite $(a_n)$ de nombres entiers (positifs à partir du rang $n=1$) et une suite $(x_n)$ de nombres réels (positifs à partir du rang $n=1$) définies par $\begin{cases} x_0:=x \\ a_0:=\entiere{x} \end{cases}$ et, pour tout $n\in\N$, 
\begin{itemize} 
\item[$\bullet$] si $x_n \in\N$, le processus s'arrête,
\item[$\bullet$] sinon, on pose $\begin{cases} x_{n+1}:=\dfrac{1}{x_n-a_n} \\ a_{n+1}:=\entiere{x_{n+1}} \end{cases}$.
\end{itemize}
Le rationnel $R_n:=\left[a_0, a_1, ..., a_n\right]$ est alors appelé la réduite d'indice $n$ dans le développement en fraction continuée du réel $x$ (ou plus simplement la réduite d'indice $n$ de $x$) et les nombres $a_0, a_1, ..., a_n$ sont appelés les coefficients de $R_n$. 

Le réel $x_n$ est appelé le quotient complet d'indice $n$ dans le développement en fraction continuée de $x$ (ou plus simplement le quotient complet d'indice $n$ de $x$) et le nombre $a_n$ est appelé le quotient partiel d'indice $n$ dans le développement en fraction continuée de $x$ (ou plus simplement le quotient partiel d'indice $n$ de $x$). 
\end{definition}

Par une récurrence immédiate, on montre que, pour tout $n\in\N$, si $x_n$ existe, alors 
\begin{equation}
x=\left[a_0, a_1, ..., a_{n-1}, x_n\right]
\label{eq070615}
\end{equation}

De plus, en conservant les notations précédentes, on montre facilement que la suite $(x_n)$ est finie si et seulement si $x\in\Q$ (\cite[p. 34]{Du07}) et, avec un peu plus de travail, que, si $x$ est un nombre irrationnel, alors la suite $(R_n)$ converge vers $x$ (\cite[p. 92-93]{Hi08}). Dans ce dernier cas, on note $x=[a_0, a_1, ..., a_n, ...]$ pour traduire l'égalité $x=\lim\limits_{n\to+\infty} R_n$. Si $x$ est irrationel, cette écriture est unique en ce sens que si $(b_n)_{n\in\N}$ est une suite d'entiers relatif strictement positif à partir du rang 1 telle que $x=[b_0, b_1, ..., b_n, ...]$ alors $a_n=b_n$ pour tout $n\in\N$ (\cite[p. 179]{HW60}, \cite[p. 36]{Du07}). On dit alors que $[a_0, a_1, ..., a_n, ...]$ est le développement de $x$ en fraction continuée. De plus, les réduites de $x$ sont, pour tout $n\in\N$, de la forme $R_n=\frac{p_n}{q_n}$ où les suites $(p_n)$ et $(q_n)$ sont définies par les relations (\ref{eq29101445}). Cette écriture est alors la forme irréductible de $R_n$ et les suites $(p_n)$ et $(q_n)$ sont également définies pour $n=-2$ et $n=-1$. Ces deux points sont considérés dans toute la suite comme implicites lorsqu'on écrit \og $R_n=\frac{p_n}{q_n}$ est la réduite d'indice $n$ de $x$ \fg{} et on ne le rappellera pas à chaque fois.

\bigskip

La propriété fondamentale des réduites dont nous nous servirons est le théorème suivant dû à A.-M. Legendre (\cite[p. 93]{Hi08}, \cite[p. 39]{Du07}).

\begin{propriete} (Legendre, 1798) --- Soit $x$ un irrationnel et $(R_n)$ la suite de ses réduites écrites sous forme irréductible $R_n=\dfrac{p_n}{q_n}$. Si un rationnel $r=\dfrac{p}{q}$ écrit sous forme irréductible vérifie $\abs{x-\dfrac{p}{q}}<\dfrac{1}{2q^2}$ alors $r=R_n$ où $n$ est l'unique entier tel que $q_{n}\leqslant q < q_{n+1}$.
\label{prop30101435}
\end{propriete}

\subsection{Nombres irrationnels quadratiques}

Il n'est en général pas simple de déterminer la forme du développement en fraction continuée d'un irrationnel. On connait celui de $\sqrt{2}=\left[1,2,2,2,...\right]$ ou celui de $e=\left[2,1,2,1,1,4,1,1,6,1,1,...,1,2n,1,...\right]$ (\cite[pp. 35-36]{Du07}) mais celui de $\pi$ reste mystérieux. La résolution de l'équation $(E_2)$ va nécessiter de déterminer le développement en fraction continuée de $\sqrt{d}$ où $d$ est un entier naturel qui n'est pas un carré parfait. Pour cela, nous allons étudier la forme générale du développement de certains nombres irrationnels: les irrationnels quadratiques réduits.

\bigskip

Un nombre irrationnel $\alpha$ est dit irrationnel quadratique s'il est racine d'un polynôme de degré 2 à coefficients entiers. On montre sans peine qu'un tel nombre s'écrit de manière unique sous la forme $\alpha=r+t\sqrt{d}$ où $r$ et $t$ sont des rationnels tels que $t\neq 0$ et $d\geqslant 2$ est un entier sans facteur carré  et que, réciproquement, tout réel de cette forme est un irrationnel quadratique.

\bigskip

\begin{definition} --- Soit $\alpha$ un irrationnel quadratique écrit sous la forme $\alpha=r+t\sqrt{d}$ avec $d\geqslant 2$ un entier sans facteur carré et $r$ et $t\neq 0$ des rationnels.
\begin{enumerate}
\item On appelle conjugué (algébrique) de $\alpha$ le nombre $\alpha^*=r-t\sqrt{d}$.
\item On dit que $\alpha$ est irrationnel quadratique réduit si $\alpha$ est un irrationnel quadratique tel que $\alpha>1$ et $-1<\alpha^*<0$.
\end{enumerate}
\end{definition}

\begin{ex} --- Le nombre d'or $\varphi:=\dfrac{1+\sqrt{5}}{2}$ est un irrationnel quadratique réduit car $\varphi>1$ et $\varphi^*=\dfrac{1-\sqrt{5}}{2}\in\intervalleoo{-1}{0}$.
\end{ex}

\bigskip

Il est clair que si $\alpha$ est un irrationnel quadratique alors $\alpha^*$ est également un irrationnel quadratique et il n'est pas difficile de voir si $P$ est un polynôme de degré 2 de $\Z[X]$ qui annule $\alpha$ alors $\alpha^*$ est l'autre racine réelle de $P$. 

En termes de structure, si $\alpha$ est un irrationnel quadratique alors $\Q[\alpha]=\enstq{u+v\alpha}{(u,v)\in\Q^2}$ est une extension de corps de $\Q$ de degré 2 i.e. un corps contenant $\Q$ et muni d'une structure de $\Q-$algèbre de dimension 2 (une base étant $(1,\alpha)$). On vérifie aisément que la conjugaison algébrique (prolongée aux rationnels en posant $u^*=u$ pour tout $u\in\Q$) est alors un automorphisme $\Q-$linéaire de corps de $\Q[\alpha]$. En particulier, pour tous rationnels $u$ et $v$, 
\begin{equation}
(u+v\alpha)^*=u+v\alpha^* \qquad \text{et} \qquad \left(u+\frac{v}{\alpha}\right)^*=u+\frac{v}{\alpha^*}. 
\label{eq07061611}
\end{equation}

Nous allons à présent caractériser les irrationnels quadratiques réduits par leur développement en fraction continuée. Commençons par rappeler les résultats obtenus par L. Euler puis par J.-L. Lagrange sur les développements en fraction continuée des irrationnels quadratiques (\cite[p. 184 et suiv.]{HW60}, \cite[pp. 40-41]{Du07}).

\begin{propriete} (Euler, 1748) --- Soit $x$ un irrationnel et $x=[a_0, a_1, ..., a_n, ...]$ son développement en fraction continuée. Si la suite $(a_n)$ est périodique à partir d'un certain rang alors $x$ est un irrationnel quadratique. 
\label{prop01110853}
\end{propriete}

\begin{propriete} (Lagrange, 1768) --- Soit $x$ un irrationnel quadratique. Alors, le développement en fraction continuée de $x$ est périodique à partir d'un certain rang.
\label{prop01110848}
\end{propriete}

\bigskip

Soit $x$ un irrationnel quadratique et $x=\left[a_0, a_1, ..., a_n, ...\right]$ son développement en fraction continuée. Si la suite $(a_n)$ est périodique de période $T$ à partir du rang $N$, on dit que le développement de $x$ est périodique de période $T$ et on note

\[x=\left[a_0, a_1, ..., a_{N-1}, \overline{a_{N}, a_{N+1}, ..., a_{N+T-1}}\right].\]

\begin{corollaire} --- Soit $x$ un irrationnel quadratique et $x=[a_0, a_1, a_2, ...]$ son développement en fraction continuée. Alors, pour tout $n\in\N$, le quotient complet $x_n=[a_n, a_{n+1}, a_{n+2}, ...]$ est un irrationnel quadratique.
Si, de plus, $x$ est réduit alors, pour tout $n\in\N$, $x_n$ est réduit et $a_n=\entiere{-\dfrac{1}{x_{n+1}^*}}$.
\label{coro18042015}
\end{corollaire}

\emph{Preuve} --- Comme le développement en fraction continuée de $x_n$ est le même que celui de $x$ aux $n$ premiers termes près, celui-ci est périodique à partir d'un certain rang donc tous les $x_n$ sont des irrationnels quadratiques d'après la propriété \ref{prop01110853}.

Supposons que $x$ soit réduit et montrons par récurrence sur $n$ que $x_n$ est également réduit. Pour $x_0=x$, c'est vrai par hypothèse. Supposons que $x_n$ soit réduit. Alors, $x_n>1$ et  $x_n^*\in\intervalleoo{-1}{0}$. \'Ecrivons $x_n=a_n +\dfrac{1}{x_{n+1}}$ (avec, rappelons-le, $a_n=\entiere{x_n}$). Alors, $x_{n+1}=\dfrac{1}{x_n-a_n}>1$ car $0<x_n-a_n<1$. De plus, d'après (\ref{eq07061611}), $x_n^*=a_n+\dfrac{1}{x_{n+1}^*}$ donc $\dfrac{1}{x_{n+1}^*}=x_n^*-a_n$. Or, $x_n>1$ donc $a_n\geqslant 1$ et $x_n^*<0$ donc $x_n^*-a_n<-1$. Ainsi, $\dfrac{1}{x_{n+1}^*}<-1$ donc $x_{n+1}^*\in\intervalleoo{-1}{0}$ ce qui montre que $x_{n+1}$ est réduit et achève la récurrence. De plus, on a démontré que, pour tout $n\in\N$, $-\dfrac{1}{x_{n+1}^*}=a_n-x_n^*$ avec $0<-x_n^*<1$ donc $a_n=\entiere{-\dfrac{1}{x_{n+1}^*}}$. \hfill $\square$

\bigskip

Nous allons pouvoir à présent donner une caractérisation des irrationnels quadratiques réduits qui est due à \'E. Galois.

\begin{propriete} (Galois, 1828) --- Un réel $x$ est irrationnel quadratique réduit si et seulement si son développement en fraction continuée $x=\left[a_0, a_1, ..., a_n, ...\right]$ est purement périodique i.e. si et seulement si la suite $(a_n)$ est périodique à partir du rang $n=0$. 
\label{prop01111315}
\end{propriete}

\emph{Preuve} --- Supposons que le développement en fraction continuée de $x$ soit purement périodique i.e. qu'il existe un entier $T\geqslant 1$ tel que $x=\left[\overline{a_0, a_1, ..., a_{T-1}}\right]$. D'après la propriété \ref{prop01110853}, $x$ est un irrationnel quadratique. De plus, $x>a_0=a_T$ et $a_T \geqslant 1$ donc $x>1$. 

Montrons que $x^*\in\intervalleoo{-1}{0}$. Notons, pour tout $n\in\N$, $\frac{p_n}{q_n}$ la réduite d'indice $n$ de $x$. \'Etant donné que $x=\left[\overline{a_0, a_1, ..., a_{T-1}}\right]$, on a également $x=[a_0, a_1, ... , a_{T-1}, x]$. Dès lors, d'après le point 3. du lemme 1,  puisque $T\geqslant 1$,
 \[x=\dfrac{p_{T-1}x+p_{T-2}}{q_{T-1}x+q_{T-2}}.\]
Ainsi, $x$ est racine du polynôme $P:=q_{T-1}X^2+(q_{T-2}-p_{T-1})X-p_{T-2}$. Comme $x$ et $x^*$ sont les deux racines de $P$ et comme $x>1$, pour montrer que $x^*\in\intervalleoo{0}{1}$, il suffit de montrer que $P$ admet une racine dans cet intervalle. Or, $P(0)=-p_{T-2}\leqslant 0$ et $P(-1)=q_{T-1}-q_{T-2}+p_{T-1}-p_{T-2}$. De plus, d'après le lemme \ref{lem29101626}, $(p_n)$ et $(q_n)$ sont croissantes à partir du rang $n=-1$ car $a_0\geqslant0$ donc, comme $T\geqslant 1$, $P(-1)\geqslant 0$. Ainsi, d'après le théorème des valeurs intermédiaires, $P$ admet une racine dans $\intervalleff{-1}{0}$. Par ailleurs, les racines de $P$ sont irrationnelles donc $P$ admet une racine dans $\intervalleoo{-1}{0}$ et donc $x^*\in\intervalleoo{0}{1}$ .

Réciproquement, supposons que $x$ soit un irrationnel quadratique réduit et notons $x=\left[a_0, a_1, ..., a_n, ...\right]$ son développement en fraction continuée. On sait, par la propriété \ref{prop01110848} que ce développement est périodique. Raisonnons par l'absurde en supposant qu'il ne soit pas purement périodique. Notons $T$ la plus petite période de ce développement et $N$ le premier indice de la première période. Ainsi, $x=\left[a_0, a_1, ..., a_{N-1}, \overline{a_N, ..., a_{N+T-1}}\right]$ avec $a_{N-1}\neq a_{N+T-1}$ et $N\geqslant 1$. Par définition de $T$, on a  $x_N=\left[a_N, a_{N+1}, ..., a_n, ...\right] = \left[a_{N+T}, a_{N+T+1}, ..., a_{n+T}, ...\right]=x_{N+T}$ donc $-\dfrac{1}{x_{N}^*}=-\dfrac{1}{x_{N+T}^*}$ et, d'après le corollaire \ref{coro18042015}, $a_{N-1}=\entiere{-\dfrac{1}{x_{N}^*}}=\entiere{-\dfrac{1}{x_{n+T}^*}}=a_{N+T-1}$ ce qui fournit la contradiction souhaitée. Ainsi, le développement de $x$ est bien purement périodique et $x=\left[\overline{a_0, a_1, ..., a_{T-1}}\right]$. \hfill $\square$

\bigskip

Nous avons vu dans le corollaire \ref{coro18042015} que les développements de $x$ et de $-\frac{1}{x^*}$ sont liés. Nous allons à présent préciser ce lien. 

\begin{propriete} --- Soit $x$ un irrationnel quadratique réduit et $\left[\overline{a_0, a_1, ..., a_{T-1}}\right]$ son développement en fraction continuée où $T$ est la plus petite période de ce développement. Alors, $-\frac{1}{x^*}$ est un irrationnel quadratique réduit et son  développement en fraction continuée est $\left[\overline{a_{T-1}, a_{T-2}, ..., a_0}\right]$.
\label{prop07061801}
\end{propriete}

\emph{Preuve} --- Posons $y=-\frac{1}{x^*}$ et notons, pour tout $n\in\N$, $y_n$ le quotient complet d'indice $n$ de $y$, $b_n$ le quotient partiel d'indice $n$ de $y$ et $x_n$ le quotient complet d'indice $n$ de $x$. On prolonge les suites $(a_n)$ et $(x_n)$ à $\Z$ en posant, pour tout $n\in\Z\setminus\N$, $a_n=a_r$ et $x_n=x_r$ où $r$ est le reste dans la division euclidienne de $n$ par $T$. Ainsi, les suites  $(a_n)_{n\in\Z}$ et $(x_n)_{n\in\Z}$ sont encore $T-$périodiques et vérifient, pour tout $n\in\Z$, $x_{n}=a_n+\frac{1}{x_{n+1}}$ (par définition) et $a_{n}=\entiere{-\frac{1}{x_{n+1}^*}}$ (d'après le corollaire \ref{coro18042015}). Montrons par récurrence sur $n$ que, pour tout $n\in\N$, $b_n=a_{T-n-1}$ et $y_{n+1}=-\frac{1}{x_{T-n-1}^*}$. \'Etant donné que $(x_n)$ est $T-$périodique, $x_T=x_0=x$ donc, puisque $x_{T-1}=a_{T-1}+\frac{1}{x_T}$, $-\frac{1}{x}=a_{T-1}-x_{T-1}$. On déduit de (\ref{eq07061611}) que $y=-\frac{1}{x^*}=\left(-\frac{1}{x}\right)^*=a_{T-1}-x_{T-1}^*$. Or, par le corollaire \ref{coro18042015}, $x_{T-1}$ est réduit donc $-x_{T-1}^*\in\intervalleoo{0}{1}$. Ainsi, $b_0=\entiere{y}=a_{T-1}$ et $y_1=\frac{1}{y-b_0}=-\frac{1}{x_{T-1}^*}$ donc le résultat est établi pour $n=0$. Supposons qu'il soit vrai pour un certain $n\in\N$. Alors, d'une part, $b_{n+1}=\entiere{y_{n+1}}=\entiere{-\frac{1}{x_{T-n-1}^*}}=a_{T-n-2}$. D'autre part, puisque $x_{T-n-2}=a_{T-n-2}+\frac{1}{x_{T-n-1}}$, 
\[y_{n+1}=-\frac{1}{x_{T-n-1}^*}=\left(-\frac{1}{x_{T-n-1}}\right)^*=a_{T-n-2}-x_{T-n-2}^*=b_{n+1}-x_{T-n-2}^*\]
donc $y_{n+2}=\frac{1}{y_{n+1}-b_{n+1}}=-\frac{1}{x_{T-n-2}^*}$ ce qui achève la récurrence.

On a donc montré que, pour tout $n\in\N$, $b_n=a_{T-n-1}$ ce qui assure que la suite $(b_n)$ est $T-$périodique donc $y$ est irrationnel quadratique réduit et, de plus, $y=\left[\overline{a_{T-1}, a_{T-2}, ..., a_0}\right]$.  \hfill $\square$

\bigskip

Cette propriété des irrationnels quadratiques réduits va servir pour étudier le développement en fraction continuée de $\sqrt{d}$ où $d$ est un entier qui n'est pas un carré parfait en fournissant un algorithme permettant de déterminer la période et de calculer les réduites de ce développement.

\begin{propriete} --- Soit $d$ un entier naturel qui n'est pas un carré parfait. On pose $x:=\entiere{\sqrt{d}}+\sqrt{d}$, on note $(x_n)$ la suite des quotients complets de $x$ et, pour tout $n\in\N$, $R_n=\frac{p_n}{q_n}$ la réduite d'indice $n$ dans le développement en fraction continuée de $\sqrt{d}$.
\begin{enumerate}
\item Le développement en fraction continuée de $\sqrt{d}$ est périodique à partir du rang 1. Plus précisément, ce développement est de la forme $\sqrt{d}=\left[a_0, \overline{a_1, a_2, ..., a_{T-1}, 2a_0}\right]$.

De plus, pour tout $k\in\llbracket 1, T-1 \rrbracket$, $a_{k}=a_{T-k}$. (Autrement dit, la liste $(a_1, a_2, ..., a_{T-1})$ est un palindrome i.e. $(a_1, a_2, ..., a_{T-1})=(a_{T-1}, a_{T-2}, ... , a_1)$).
\item Il existe deux suites $T-$périodiques d'entiers naturels non nuls $(u_n)$ et $(v_n)$ telles que, pour tout entier $n\in\N$, $x_n=\dfrac{u_n+\sqrt{d}}{v_n}$, $u_n\leqslant a_0$ et $v_n$ divise $d-u_n^2$.
\item Pour tout $n\geqslant 0$, $p_{n-1}^2-dq_{n-1}^2=(-1)^nv_n$.
\item La plus petite période $T$ dans le développement de $\sqrt{d}$ est le premier indice $k$ tel que $a_k=2a_0$.
\item S'il existe un entier $N\in\llbracket 1, T-1\rrbracket$ tel que $v_N$ divise $u_N$ alors $T=2N$ et, en particulier, $T$ est pair.
\end{enumerate}
\label{prop01112157}
\end{propriete}

\emph{Preuve}
\begin{enumerate}
\item Notons $\sqrt{d}=\left[a_0, a_1, ..., a_n, ...\right]$ le développement en fraction continuée de $\sqrt{d}$. On a donc, en particulier, $a_0=\entiere{\sqrt{d}}$. \'Ecrivons, de même, $x=\left[a_0', a_1', ..., a_n', ...\right]$ le développement de $x$. Remarquons que, comme $d$ n'est pas un carré parfait, $\sqrt{d}$ est un irrationnel quadratique et donc $x$ aussi. De plus, comme $d$ n'est pas un carré, $d\geqslant 2$ donc $x>1$ et $x^*=a_0-\sqrt{d} \in\intervalleoo{-1}{0}$ puisque $a_0=\entiere{\sqrt{d}}$. Ainsi, $x$ est réduit. On déduit alors de la propriété \ref{prop01111315} que son développement est purement périodique. De plus, comme $a_0'=\entiere{x}=2a_0$, ce développement est de la forme $x=\left[\overline{2a_0, a_1', ..., a_{T-1}'}\right]$ en notant $T$  la plus petite période de ce développement. On en déduit que le développement de $\sqrt{d}=x-a_0$ est $\left[a_0, \overline{a_1', a_2', ..., a_{T-1}', 2a_0}\right]=\left[a_0, \overline{a_1, a_2, ..., a_{T-1}, 2a_0}\right]$ car, par construction de $x$, pour tout $n\geqslant 1$, $a_n'=a_n$.

Par définition, $x_1=\frac{1}{x-\entiere{x}}=\frac{1}{\sqrt{d}-a_0}$ donc $-\frac{1}{x_1^*}=\left(-\frac{1}{x_1}\right)^*=\sqrt{d}+a_0=x$. Or, $x_1=\left[\overline{a_1', a_2', ..., a_{T_1}', 2a_0}\right]$ donc, $x_1$ étant un irrationnel quadratique réduit d'après le corollaire \ref{coro18042015}, on déduit de la propriété \ref{prop07061801} que  $-\frac{1}{x_1^*}=\left[\overline{2a_0, a_{T-1}', a_{T-2}', ..., a_{1}'}\right]$. Ainsi, $\left[\overline{2a_0, a_1', ..., a_{T-1}'}\right]=\left[\overline{2a_0, a_{T-1}', a_{T-2}', ..., a_{1}'}\right]$ et donc, pour tout $k\in\llbracket 1, T-1 \rrbracket$, $a_k'=a_{T-k}'$. \'Etant donné que $a_n=a_n'$ pour tout $n\geqslant 1$, on conclut que $(a_1, a_2, ..., a_{T-1})$ est un palindrome.

\item Construisons les suites $(u_n)$ et $(v_n)$ par récurrence. Pour $n=0$, $x_0=x=a_0+\sqrt{d}$ donc $u_0=a_0$ et $v_0=1$. Supposons qu'on a construit $u_n$ et $v_n$ pour un certain $n\in\N$. Alors, en rappelant que $a_n'=\entiere{x_n}$,
\[x_{n}=\dfrac{u_n+\sqrt{d}}{v_n}=a_n'+\dfrac{u_n-a_n'v_n+\sqrt{d}}{v_n} = a_n'+\dfrac{-u_{n+1}+\sqrt{d}}{v_n}\]
où on a posé $u_{n+1}:=a_n'v_n-u_n$. Dès lors, 
\[x_{n}=a_n'+\dfrac{d-u_{n+1}^2}{v_n(u_{n+1}+\sqrt{d})}=a_n'+\dfrac{1}{\dfrac{u_{n+1}+\sqrt{d}}{\frac{d-u_{n+1}^2}{v_n}}}=a_n'+\dfrac{1}{\dfrac{u_{n+1}+\sqrt{d}}{v_{n+1}}}\]
en posant $v_{n+1}:=\dfrac{d-u_{n+1}^2}{v_n}$. Ainsi, on peut affirmer que $x_{n+1}=\dfrac{u_{n+1}+\sqrt{d}}{v_{n+1}}$. Par définition, $u_{n+1}$ est un entier. De plus, comme $u_n\leqslant a_0<\sqrt{d}$, $x_n=\dfrac{u_n+\sqrt{d}}{v_n}>\dfrac{2u_n}{v_n}$. Or, d'après le corollaire \ref{coro18042015}, $x_n$ est réduit donc $x_n>1$ et ainsi $a_n'=\entiere{x_n}>\dfrac{x_n}{2}>\dfrac{u_n}{v_n}$ ce qui assure que $u_{n+1}>0$. Par ailleurs, $u_{n+1}=a_n'v_n-u_n<x_nv_n-u_n=\sqrt{d}$ donc $u_{n+1}\leqslant a_0$ et $v_{n+1}>0$. Il ne reste plus qu'à montrer que $v_n$ divise $d-u_{n+1}^2$ ce qui prouvera à la fois que $v_{n+1}$ est entier et que $v_{n+1}$ divise $d-u_{n+1}^2$. Pour le voir, il suffit d'écrire $d-u_{n+1}^2=d-(a_n'v_n-u_n)^2=d-u_n^2-v_n({a_n'}^2v_n-2a_n'u_n)$ et d'utiliser le fait que, par hypothèse de récurrence, $v_n$ divise $d-u_n^2$.

Ainsi, les suites $(u_n)$ et $(v_n)$ sont définies par $u_0=a_0=\entiere{\sqrt{d}}$, $v_0=1$ et, pour tout $n\in\N$,
\begin{equation}
u_{n+1}=a_n'v_n-u_n=\begin{cases} a_0 \text{ si }n=0 \\ a_nv_n-u_n \text{ si }n\geqslant 1 \end{cases} \text{ et }~~ v_{n+1}=\frac{d-u_{n+1}^2}{v_n}
\label{eq20041953}
\end{equation}

Pour finir, notons que $(x_n)$ est $T-$périodique car le développement de $x$ l'est. Ainsi, pour tout $n\in\N$, $x_{n+T}=x_n$ i.e. $\frac{u_{n+T}}{v_{n+T}}+\frac{1}{v_{n+T}}\sqrt{d}=\frac{u_{n}}{v_{n}}+\frac{1}{v_{n}}\sqrt{d}$ donc, comme $(1,\sqrt{d})$ est une $\Q-$base de $\Q[\sqrt{d}]$, $v_{n+T}=v_n$ et $u_{n+T}=u_n$ et ainsi les suites $(u_n)$ et $(v_n)$ sont $T-$périodiques.
\item Soit $n\in\N$. Si $n=0$ alors $p_{n-1}^2-dq_{n-1}^2=1^2-d\times 0^2=1=v_0$.

Si $n\geqslant 1$ alors, le développement en fraction continuée de $\sqrt{d}$ est
\[\sqrt{d}=\left[a_0, a_1, ..., a_{n-1}, x_n\right]\]
donc, d'après le lemme \ref{lem29101626},
\[\sqrt{d}=\dfrac{p_{n-1}x_n+p_{n-2}}{q_{n-1}x_n+q_{n-2}}=\dfrac{p_{n-1}(u_n+\sqrt{d})+v_np_{n-2}}{q_{n-1}(u_n+\sqrt{d})+v_nq_{n-2}}.\]
Il s'ensuit que
\[(q_{n-1}u_n+v_nq_{n-2})\sqrt{d}+dq_{n-1}=p_{n-1}\sqrt{d}+p_{n-1}u_n+v_np_{n-2}\]
donc, en utilisant encore le fait que $(1,\sqrt{d})$ est une $\Q-$base de $\Q[\sqrt{d}]$,
\begin{equation} \begin{cases} p_{n-1}=q_{n-1}u_n+q_{n-2}v_n \\ dq_{n-1}=p_{n-1}u_n+p_{n-2}v_n \end{cases}.
\label{eq02111252}
\end{equation}
On en déduit, en utilisant également le point 2. du lemme \ref{lem29101626}, que
\begin{align*}
p_{n-1}^2-dq_{n-1}^2&=p_{n-1}(q_{n-1}u_n+q_{n-2}v_n)-q_{n-1}(p_{n-1}u_n+p_{n-2}v_n) \\
&=(p_{n-1}q_{n-2}-q_{n-1}p_{n-2})v_n=(-1)^nv_n.
\end{align*}
Ainsi, la relation est vraie pour tout $n\geqslant 0$. 
\item Montrons que $T$ est le premier indice tel que $a_{T}=2a_0$. Il est équivalent de montrer qu'aucun des coefficients $a_1$, $a_2$, ..., $a_{T-1}$ n'est égal à $2a_0$. Pour cela, remarquons tout d'abord que, pour tout $j\in\llbracket 1, T-1 \rrbracket$, $x_j\neq x$. En effet, dans le cas contraire, on aurait $x=\left[2a_0, a_1, ..., a_{j-1}, x\right]=\left[\overline{2a_0, a_1, ..., a_{j-1}}\right]$ et donc le développement serait périodique de période $j<T$ ce qui contredit la définition de $T$. Supposons alors que $j\geqslant 1$ soit un indice quelconque tel que $a_j=2a_0$. Ainsi, comme $j\geqslant 1$, d'après (\ref{eq20041953}), $u_{j+1}=2a_0v_j-u_j$. Or, $u_j\leqslant a_0$ et $u_{j+1}\leqslant a_0$ donc, comme $v_j\in\N^*$, $v_j=1$ et $u_{j+1}=u_j=a_0$. On conclut donc que $x_j=a_0+\sqrt{d}=x$ et ainsi, d'après la remarque précédente, $j\notin\llbracket 1, T-1 \rrbracket$ ce qui achève la démonstration.
\item Supposons que $N\in\llbracket 1, T \rrbracket$ soit tel que $v_N$ divise $u_N$. On va montrer par une récurrence finie que, pour tout $j\in\llbracket 0, N-1\rrbracket$, on a les égalités
\begin{equation}
u_{N+j+1}=u_{N-j}, \qquad v_{N+j+1}=v_{N-j-1} \qquad \text{et} \qquad a_{N+j+1}'=a_{N-j-1}'
\label{eq20042001}
\end{equation}
Pour $j=0$, sachant que $v_N$ divise $u_N$, $\dfrac{2u_N}{v_N}$ est entier. Or, $u_N<\sqrt{d}$ donc $\dfrac{2u_N}{v_N}<\dfrac{u_N+\sqrt{d}}{v_N}$ et, d'après le corollaire \ref{coro18042015}, $x_N=\dfrac{u_N+\sqrt{d}}{v_N}$ est un irrationnel quadratique réduit donc $x_N^*\in\intervalleoo{-1}{0}$ i.e. $\dfrac{u_N-\sqrt{d}}{v_N}\in\intervalleoo{-1}{0}$. Dès lors, $\dfrac{u_N+\sqrt{d}}{v_N}=\dfrac{2u_N}{v_N}-\dfrac{u_N-\sqrt{d}}{v_N}<\dfrac{2u_N}{v_N}+1$. Il s'ensuit que $a_n'=\entiere{\dfrac{u_N+\sqrt{d}}{v_N}}=\dfrac{2u_N}{v_N}$. Or, d'après (\ref{eq20041953}), $u_{N+1}=a_N'v_N-u_N$ donc $u_{N+1}=u_{N}$. De plus, sachant que $u_{N+1}^2+v_{N+1}v_N=d=u_N^2+v_Nv_{N-1}$, on peut affirmer que $v_{N+1}=v_{N-1}$. Dès lors,
\[a_{N+1}'=\entiere{x_{N+1}}=\entiere{\dfrac{u_{N+1}+\sqrt{d}}{v_{N+1}}}=\entiere{\dfrac{u_{N}+\sqrt{d}}{v_{N-1}}}.\]
Or, $u_{N}=a_{N-1}'v_{N-1}-u_{N-1}$ donc 
\[a_{N+1}'=\entiere{a_{N-1}'-\dfrac{u_{N-1}-\sqrt{d}}{v_{N-1}}}=\entiere{a_{N-1}'-x_{N-1}^*}=a_{N-1}'\]
car $x_{N-1}$ est réduit donc $x_{N-1}^*\in\intervalleoo{-1}{0}$.

Supposons que les relations (\ref{eq20042001}) soient vraies pour un certain $j\in\llbracket 0, N-2\rrbracket$. Alors, 
\begin{align*}
u_{N+j+2}&=a_{N+j+1}'v_{N+j+1}-u_{N+j+1}=a_{N-j-1}'v_{N-j-1}-u_{N-j} \\ 
&=a_{N-j-1}'v_{N-j-1}-(a_{N-j-1}'v_{N-j-1}-u_{N-j-1})\end{align*}
donc $u_{N+j+2}=u_{N-j-1}$. Par suite, en écrivant que 
\[u_{N+j+2}^2-v_{N+j+2}v_{N+j+1}=d=u_{N-j-1}^2+v_{N-j-1}v_{N-j-2}=u_{N+j+2}^2+v_{N+j+1}v_{N-j-2},\] on est assuré que $v_{N+j+2}=v_{N-j-2}$. Enfin, en raisonnant comme précédemment,
\begin{align*}
a_{N+j+2}'&=\entiere{x_{N+j+2}}=\entiere{\dfrac{u_{N+j+2}+\sqrt{d}}{v_{N+j+2}}}=\entiere{\dfrac{u_{N-j-1}+\sqrt{d}}{v_{N-j-2}}}\\
&=\entiere{a_{N-j-2}'-\dfrac{u_{N-j-2}-\sqrt{d}}{v_{N-j-2}}}=\entiere{a_{N-j-2}'-x_{N-j-2}^*}\end{align*}
et donc $a_{N+j+2}=a_{N-j-2}'$ car $x_{N-j-2}^*\in\intervalleoo{-1}{0}$ ce qui achève la récurrence.

Supposons  à présent que $N\in\llbracket 1, T-1 \rrbracket$. Si on applique ce qui précède avec $j=N-1$, on obtient, en particulier, que $a_{0}'=a_{2N}'$ i.e. comme $N\geqslant 1$, $2a_0=a_{2N}$. Ainsi, d'après le point 3., sachant que $2\leqslant 2N \leqslant 2(T-1) < 2T$, on en déduit que $2N=T$.
 \hfill  $\square$
\end{enumerate}

\bigskip

\begin{ex} --- Pour un entier donné $d$ qui n'est pas un carré parfait, la démonstration nous donne non seulement un critère pour déterminer une période minimale dans le développement en fraction continuée de $\sqrt{d}$ (on s'arrête dès qu'un coefficient d'une réduite est égale à $2\entiere{d}$) mais elle fournit, de plus, un algorithme pour déterminer ce développement. Par exemple, si $d=21$, on part de $x=\entiere{\sqrt{21}}+\sqrt{21}=4+\sqrt{21}$ puis on écrit successivement:

\begin{align*}
x&=8+(-4+\sqrt{21})=8+\dfrac{5}{4+\sqrt{21}}=\boxed{8}+\dfrac{1}{x_1} \\
x_1&=\dfrac{4+\sqrt{21}}{5}=1+\dfrac{-1+\sqrt{21}}{5}=1+\dfrac{20}{5(1+\sqrt{21})}=\boxed{1}+\dfrac{1}{x_2} \\
x_2&=\dfrac{1+\sqrt{21}}{4}=1+\dfrac{-3+\sqrt{21}}{4}=1+\dfrac{12}{4(3+\sqrt{21})}=\boxed{1}+\dfrac{1}{x_3} \\
x_3&=\dfrac{3+\sqrt{21}}{3}=2+\dfrac{-3+\sqrt{21}}{3}=2+\dfrac{12}{3(3+\sqrt{21})}=\boxed{2}+\dfrac{1}{x_4} \\
x_4&=\dfrac{3+\sqrt{21}}{4}=1+\dfrac{-1+\sqrt{21}}{4}=1+\dfrac{20}{4(1+\sqrt{21})}=\boxed{1}+\dfrac{1}{x_5} \\
x_5&=\dfrac{1+\sqrt{21}}{5}=1+\dfrac{-4+\sqrt{21}}{5}=1+\dfrac{5}{5(4+\sqrt{21})}=\boxed{1}+\dfrac{1}{x}
\end{align*}
donc $x=\left[8,1,1,2,1,1,x\right]=\left[\overline{8,1,1,2,1,1}\right]$ et ainsi $\sqrt{21}=\left[4,\overline{1,1,2,1,1,8}\right]$. 
\label{ex01112333}
\end{ex}

\section{L'équation diophantienne $\boldsymbol{ax^2-by^2=1}$}
\subsection{Un cas particulier: les équations de Pell-Fermat}
\label{secPell}
On appelle \emph{équation de Pell-Fermat} une équation diophantienne d'inconnue $(x,y)\in\N^2$ de la forme
\[ (F_d): x^2-dy^2=1\]
où $d$ est un entier naturel non nul.

L'étude d'une telle équation remonte à l'antiquité (Voir l'article \cite{Le02} qui aborde le problème dit des \emph{boeufs d'Hélios} qui aurait été posé par Archimède à \'Eratosthène). Elle n'a cependant trouvé de résolution définitive qu'au XVIIIe siècle. Il semble que le nom de Pell a été attaché à cette équation suite à une erreur de L. Euler qui a attribué à J. Pell une méthode lue dans un ouvrage d'un autre mathématicien anglais, en l'occurrence J. Wallis. En Europe, la première méthode de résolution est due à Lord W. Brouncker en 1657 mais les mathématiciens indiens Brahmagupta et Bh\={a}skara II (aussi appelé Bh\={a}skar\={a}c\={a}rya) avaient mis au point, entre le VIIe et le XIIe siècle, un algorithme de résolution connu aujourd'hui sous le nom de \emph{méthode chakravala}. Par la suite, J. Wallis et P. de Fermat ont été les premiers à affirmer que l'équation $(F_d)$ avait toujours une infinité de solutions lorsque $d$ n'est pas le carré d'un entier avant que la démonstration ne soit donnée par J.-L. Lagrange en 1766. Pour une discussion exhaustive sur l'histoire de l'équation de Pell-Fermat, voir \cite{Whi12}.

Pour tout entier $d$, l'équation $(F_d)$ admet pour solution évidente $(x, y)=(1, 0)$. Nous dirons que c'est une solution triviale de $ (F_d)$. Une façon rapide de démontrer que l'équation $(F_d)$ admet toujours des solutions non triviales dans le cas où  $d$ est un entier naturel qui n'est pas un carré parfait est d'utiliser le théorème de Dirichlet (\cite[p. 89]{Hi08}, \cite[p. 5 et p. 38]{Du07}). Cette méthode a cependant l'inconvénient d'être non effective: elle ne donne pas de moyen de déterminer les solutions. On peut compléter ce point de vue en montrant l'existence d'une solution \og fondamentale \fg{} $(x_1, y_1)$ de $(F_d)$ telle que les solutions de $(F_d)$ sont exactement les couples $(x_n, y_n)$ définis par $\left(x_1+y_1\sqrt{d}\right)^n=x_n+y_n\sqrt{d}$ (\cite[p. 53-54]{Mo69} ou \cite[p. 264]{Hu82}), propriété intimement liée à la théorie des unités des corps quadratiques (\cite[p. 51]{Du07}). Cette méthode ne permet cependant de conclure que dans le cas où on sait déterminer la solution fondamentale.

Nous proposons ici une méthode effective de résolution de $(F_d)$ qui a, de plus, l'avantage de s'appliquer à des équations plus générales. 

\bigskip

On appelle \emph{équation de Pell-Fermat généralisée} une équation d'inconnue $(x,y)\in\N^2$ de la forme
\[(F_{d,m}): x^2-dy^2=m\]
où $d$ et $m$ sont des entiers relatifs non nuls et premiers entre eux.

Si $d$ est négatif alors $F_{d,m}$ n'a qu'un nombre fini (éventuellement nul) de solutions qu'on peut déterminer, par exemple, par une méthode de crible. Si $d$ est un carré parfait i.e. s'il existe un entier $k$ tel que $d=k^2$ alors $(F_{d,m})$ équivaut à $(x-ky)(x+ky)=m$ donc la résolution se résume à chercher les diviseurs de $m$. Ces cas n'ayant pas beaucoup d'intérêt, on suppose que $d$ est un entier naturel qui n'est pas un carré parfait. Il est clair que si $(x, y)$ est une solution de $(F_{d,m})$ telle que $x$ et $y$  ne sont pas premiers entre eux, alors le carré de $\delta:=\text{PGCD}(x,y)$ divise $m$ donc en divisant $(F_{d,m})$ par $\delta^2$, on est ramené à résoudre une certaine équation $(F_{d,m'})$ avec $m'<m$ dont les inconnues sont des entiers premiers entre eux. On peut donc toujours supposer que $\delta=1$.

Remarquons que, même sous ces conditions,  l'équation $(F_{d,m})$ n'a pas toujours de solution. Par exemple, pour $(d,m)=(7,5)$, il ne peut y avoir de solution car on aurait alors $x^2\equiv 5~[7]$ et les restes possibles d'un carré modulo 7 sont 0, 1, 2 ou 4. 

Dans ce qui suit, nous allons donner une condition nécessaire et suffisante pour que $(F_{d,m})$ ait une solution dans le cas où $1\leqslant\abs{m}<\sqrt{d}$ et nous montrerons qu'alors $(F_{d,m})$ a une infinité de solutions. On peut toujours imposer que $x$ et $y$ soient non nuls. Le cas $x=0$ ne peut pas se produire car alors $-dy^2=m$ ce qui est absurde car $1\leqslant \abs{m} < \sqrt{d}$ et le cas $y=0$ ne se produit que si  $m$ est un carré parfait et, dans ce cas, on dira que la solution $(\sqrt{m}, 0)$ est une solution triviale de $(F_{d,m})$.

\bigskip

\begin{lemme} --- Soit $d$ un entier naturel qui n'est pas un carré parfait et $m$ un entier relatif tel que $1\leqslant \abs{m} <\sqrt{d}$. Si $x$ et $y$ sont deux entiers naturels non nuls et premiers entre eux tels que $(x,y)$ est une solution de $(F_{d,m})$ alors $\dfrac{x}{y}$ est une réduite dans le développement en fraction continuée de $\sqrt{d}$.
\label{lem01112220}
\end{lemme}

\emph{Preuve} --- Soit $x$ et $y$ deux entiers naturels non nuls premiers entre eux tels que $x^2-dy^2=m$. Alors, $x-y\sqrt{d}=\dfrac{m}{x+y\sqrt{d}}$.

Si $m>0$ alors $x^2=m+dy^2>dy^2$ donc $x>y\sqrt{d}$ et ainsi 
\[\abs{\sqrt{d}-\dfrac{x}{y}}=\dfrac{\abs{m}}{xy+y^2\sqrt{d}}<\dfrac{\sqrt{d}}{xy+y^2\sqrt{d}} <\dfrac{\sqrt{d}}{2y^2\sqrt{d}}=\dfrac{1}{2y^2}.\]
Il s'ensuit, d'après la propriété \ref{prop30101435}, que $\dfrac{x}{y}$ est une réduite dans le développement de $\sqrt{d}$ en fraction continuée (car $\sqrt{d}$ est irrationnel).

Si $m<0$ alors $dy^2=x^2-m>x^2$ donc $y\sqrt{d}>x$ et ainsi
\[\abs{\dfrac{1}{\sqrt{d}}-\dfrac{y}{x}}=\dfrac{\abs{m}}{\sqrt{d}(x^2+xy\sqrt{d})}<\dfrac{\sqrt{d}}{\sqrt{d}(2x^2)} =\dfrac{1}{2x^2}.\]
De même, $\dfrac{y}{x}$ est une réduite dans le développement en fraction continuée de $\sqrt{d}^{\phantom{.}-1}$. \'Ecrivons alors $\sqrt{d}=\left[a_0, a_1, ... , a_n, ...\right]$ et $\sqrt{d}^{\phantom{.}-1}=\left[a_0', a_1', ... , a_n', ...\right]$ les développements en fraction continuée de $\sqrt{d}$ et $\sqrt{d}^{\phantom{.}-1}$ et notons  $(\delta_n')$ la suite des quotients complets de $\sqrt{d}^{\phantom{.}-1}$. Comme $\sqrt{d}>1$, $a_0'=\entiere{\sqrt{d}^{\phantom{.}-1}}=0$ donc $\delta_1'=0+\dfrac{1}{\sqrt{d}^{\phantom{.}-1}-0}=\sqrt{d}$ et ainsi $a_1'=\entiere{\delta_1'}=\entiere{\sqrt{d}}=a_0$. Il s'ensuit par récurrence que, pour tout $n\geqslant 1$, $a_n'=a_{n-1}$ et ainsi le développement en fraction continuée de $\sqrt{d}^{\phantom{.}-1}$ est $\left[0,a_0, a_1, ..., a_{n-1}, ...\right]$. Si on note $k$ l'entier tel que $\dfrac{y}{x}$ soit la réduite d'indice $k$ de ${\sqrt{d}}^{\phantom{.}-1}$ alors $k\neq0$ car $y\neq0$. Ainsi, 
\[\dfrac{y}{x}=\left[0,a_0, a_1, ..., a_{k-1}\right]=0+\dfrac{1}{\left[a_0, a_1, ..., a_{k-1}\right]}\]
donc $\dfrac{x}{y}=\left[a_0, a_1, ..., a_{k-1}\right]$ i.e. $\dfrac{x}{y}$ est la réduite d'indice $k-1$ dans le développement en fraction continuée de $\sqrt{d}$. \hfill $\square$

\bigskip

Cette propriété d'approximation diophantienne étant établie, nous allons pouvoir décrire complètement l'ensemble des solutions de $(F_{d,m})$ lorsque $1\leqslant\abs{m}<\sqrt{d}$.

\bigskip

\begin{theoreme} --- Soit $d$ un entier naturel qui n'est pas un carré parfait et soit $m$ un entier relatif tel que $1\leqslant \abs{m} < \sqrt{d}$. On considère la suite $(v_n)$ associée à $\sqrt{d}$ par la propriété \ref{prop01112157}, $T$ la période minimale du développement de $\sqrt{d}$ en fraction continuée et $\dfrac{p_n}{q_n}$ la réduite d'indice $n$ de ce développement. On note $\mathcal{S}_{d,m}$ l'ensemble des solutions non triviales $(x, y)$ de $(F_{d,m})$ telles que $\mathrm{PGCD}(x,y)=1$ et, pour tout $\ell\in\Z$, $\mathscr{E}_{\ell}:=\enstq{j\in\llbracket 1, T \rrbracket}{(-1)^jv_j=\ell}$. 
\begin{enumerate}
\item Si $T$ est pair, $\mathcal{S}_{d,m}$ est non vide si et seulement si $\mathscr{E}_m$ est non vide et, dans ce cas, $\mathcal{S}_{d,m}$ est la réunion des ensembles $\enstq{(p_{N+kT-1}, q_{N+kT-1})}{k\in\N}$ pour $N$ parcourant $\mathscr{E}_m$.
\item Si $T$ est impair, $\mathcal{S}_{d,m}$ est non vide si et seulement si l'un, au moins, des deux ensembles $\mathscr{E}_m$ ou $\mathscr{E}_{-m}$ est non vide et, dans ce cas, $\mathcal{S}_{d,m}$ est la réunion des ensembles $\enstq{(p_{N+2kT-1}, q_{N+2kT-1})}{k\in\N}$ pour $N$ parcourant $\mathscr{E}_m$ et des ensembles $\enstq{(p_{M+(2k+1)T-1}, q_{M+(2k+1)T-1})}{k\in\N}$ pour $M$ parcourant $\mathscr{E}_{-m}$.
\item Dans tous les cas, si $\mathcal{S}_{d, m}$ est non vide alors $(F_{d,m})$ admet une infinité de solutions.
\end{enumerate}
\label{theo01112218}
\end{theoreme}

\emph{Preuve}
\begin{enumerate}
\item On se place dans le cas où $T$ est pair. Supposons que $(x,y)\in\mathcal{S}_{d,m}$. On sait par le lemme \ref{lem01112220} que $\dfrac{x}{y}$ est une réduite de $\sqrt{d}$. Ainsi, il existe un entier $j\in\N^*$ tel que $x=p_{j-1}$ et $y=q_{j-1}$. Il s'ensuit, d'après le point 3. de la propriété \ref{prop01112157}, que $m=x^2-dy^2=p_{j-1}^2-dq_{j-1}^2=(-1)^{j}v_{j}$. De plus, si on note $k$ le quotient dans la division euclidienne de $j-1$ par $T$, il existe un entier $N_j\in\llbracket 1, T\rrbracket$ tel que $j=N_j+kT$. Comme $(v_n)$ est $T-$périodique et comme $T$ est un nombre pair, 
\[(-1)^{N_j}v_{N_j}=(-1)^{N_j+kT}v_{N_j+kT}=(-1)^{j}v_{j}=m\]
et, ainsi, $N_j\in\mathscr{E}_m$ donc $\mathscr{E}_m$ n'est pas vide. De plus, $(x, y)=(p_{N_j+kT-1}, q_{N_j+kT-1})$ appartient à la réunion des ensembles $\enstq{(p_{N+kT-1}, q_{N+kT-1})}{k\in\N}$ pour $N$ parcourant $\mathscr{E}_m$.

Réciproquement, si $\mathscr{E}_m$ n'est pas vide et si $(x,y)$ est de la forme $(p_{N+kT-1}, q_{N+kT-1})$ avec $N\in\mathscr{E}_m$ et $k\in\N$ alors, comme $(v_n)$ est $T-$périodique et comme $T$ est pair, d'après le point 3. de la propriété \ref{prop01112157},
\[x^2-dy^2=p_{N+kT-1}^2-dq_{N+kT-1}^2=(-1)^{N+kT}v_{N+kT}=(-1)^{N}v_{N}=m\]
donc $(x,y)\in\mathcal{S}_{d,m}$.
\item On se place à présent dans le cas où $T$ est impair. Supposons que $(x,y)\in\mathcal{S}_{d,m}$. Alors, on montre, comme précédemment, qu'il existe un entier $j\in\N^*$ tel que $x=p_{j-1}$, $y=q_{j-1}$ et $m=(-1)^{j}v_{j}$ et qu'on peut écrire  $j=N_j+KT$ avec $K=\entiere{\frac{j-1}{T}}$ et $N_j\in\llbracket 1, T\rrbracket$. Deux cas sont alors possibles.

\emph{1\up{er} cas}: si $K=2k$ est pair alors $KT$ est pair et on peut raisonner comme dans le point 1. pour montrer que $N_j\in\mathscr{E}_m$ donc $\mathscr{E}_m$ n'est pas vide et que $(x, y)=(p_{N_j+2kT-1}, q_{N_j+2kT-1})$ appartient à la réunion des ensembles $\enstq{(p_{N+2kT-1}, q_{N+2kT-1})}{k\in\N}$ pour $N$ parcourant $\mathscr{E}_m$.

\emph{2\up{e} cas}: si $K=2k+1$ est impair alors, comme $T$ est impair, $KT$ est également impair. Dans ce cas, comme $(v_n)$ est $T-$périodique, 
\[(-1)^{N_j}v_{N_j}=-(-1)^{N_j+(2k+1)T}v_{N_j+(2k+1)T}=-(-1)^{j}v_{j}=-m\]
et, ainsi, $N_j\in\mathscr{E}_{-m}$ donc $\mathscr{E}_{-m}$ n'est pas vide. De plus, $(x, y)=(p_{N_j+(2k+1)T-1}, q_{N_j+(2k+1)T-1})$ appartient à la réunion des ensembles $\enstq{(p_{M+(2k+1)T-1}, q_{M+(2k+1)T-1})}{k\in\N}$ pour $M$ parcourant $\mathscr{E}_{-m}$.

Réciproquement, si $\mathscr{E}_m$ n'est pas vide et si $(x,y)$ est de la forme $(p_{N+2kT-1}, q_{N+2kT-1})$ avec $N\in\mathscr{E}_m$ et $k\in\N$ alors on montre comme en 1. que $(x,y)\in\mathcal{S}_{d,m}$.  Si $\mathscr{E}_{-m}$ n'est pas vide et si $(x,y)$ est de la forme $(p_{M+(2k+1)T-1}, q_{M+(2k+1)T-1})$ avec $M\in\mathscr{E}_{-m}$ et $k\in\N$ alors, comme $(v_n)$ est $T-$périodique et comme $(2k+1)T$ est impair,
\[x^2-dy^2=(-1)^{M+(2k+1)T}v_{M+(2k+1)T}=-(-1)^{M}v_{M}=-(-m)=m\]
donc $(x,y)\in\mathcal{S}_{d,m}$.
\item Le dernier point est une conséquence directe des deux précédents et du fait que les réduites sont deux à deux distinctes. \hfill $\square$
\end{enumerate} 

\bigskip

\begin{rem} --- L'énoncé du théorème précédent est lourd mais sa philosophie est assez simple: si $T$ est pair, une solution de $x^2-dy^2=m$ ne peut provenir, à un certain nombre de périodes près, que d'un indice $N\in\llbracket 1, T \rrbracket$ tel que $(-1)^Nv_N=m$. En revanche, si $T$ est impair, une solution de $x^2-dy^2=m$ peut provenir soit d'un indice $N\in\llbracket 1, T \rrbracket$ tel que $(-1)^Nv_N=m$ à un nombre pair de périodes près soit d'un indice $M\in\llbracket 1, T \rrbracket$ tel que $(-1)^Mv_M=-m$ à un nombre impair de périodes près. Ainsi, dans la pratique, il est assez simple de savoir dans quel cas on se trouve et de déterminer les premières solutions de $(F_{d,m})$ s'il y en a.

Par ailleurs, il faut garder à l'esprit que ce théorème ne traite pas le cas où les inconnues $x$ et $y$ ne sont pas des nombres premiers entre eux.
\end{rem}

\bigskip

\begin{ex} --- Considérons l'équation $(F_{21,m}): x^2-21y^2=m$ avec $1\leqslant \abs{m} \leqslant 4$. On a vu dans l'exemple \label{ex01112333} que la suite $(v_n)$ associée à $\sqrt{21}$ est périodique de période $T=6$ avec $v_1=v_5=5$, $v_2=v_4=4$, $v_3=3$ et $v_6=1$. On est donc dans le cas le plus simple ($T$ est pair) et $(F_{21,m})$ (avec $\abs{m}\leqslant4$) admet des solutions en nombres premiers entre eux si et seulement si $m=1$, $m=4$ ou $m=-3$.

Pour $m=1$ qui correspond à $N=6$, on calcule $R_5=[4,1,1,2,1,1]=\dfrac{55}{12}$ qui fournit la solution $(55,12)$. C'est la solution fondamentale de l'équation de Pell-Fermat $(F_{21})$. Les autres solutions non triviales ($(x_2, y_2)=(\np{6049}, \np{1320})$, $(x_3, y_3)=(\np{665335}, \np{145188})$, etc...) peuvent se déterminer en utilisant les autres réduites $R_{5+6n}$ ou se déduire de l'identité $x_n+y_n\sqrt{21}=(55+12\sqrt{21})^n$.  Comme $m=1$ est un carré parfait, il y également la solution triviale $(x_0, y_0)=(1, 0)$.

Pour $m=4$ qui correspond à $N=2$ ou $N=4$, on calcule $R_1=[4,1]=\dfrac{5}{1}$ et $R_3=[4,1,1,2]=\dfrac{23}{5}$ qui donnent les solutions $(5,1)$ et $(23, 5)$ et les autres solutions correspondent aux réduites $R_{7}$, $R_9$, $R_{13}$, $R_{15}$, etc... Comme $m=4$ est un carré parfait, il y a également la solution triviale $(2,0)$.

Pour $m=-3$ qui correspond à $N=3$, on calcule $R_2=[4,1,1]=\dfrac{9}{2}$ donc $(9,2)$ est solution de $(F_{21,-3})$. Les autres solutions correspondent aux réduites $R_{8}$, $R_{14}$, etc...

Si $x$ et $y$ ne sont pas premiers entre eux et si on note $\delta$ leur P.G.C.D alors $\delta^2$ divise $m$ donc cela impose $m=4$ et $\delta=2$. Il s'ensuit que $x=2x'$ et $y=2y'$ avec $(x',y')$ solution de $X^2-21Y^2=1$ i.e. de $(F_{21})$ et on est ramené au premier cas qui va donner les solutions $(110,24)$, $(\np{12098},\np{2640})$, ... pour l'équation $(F_{21,4})$.
\end{ex}

\bigskip

Pour terminer, nous en déduisons une démonstration du théorème de Lagrange concernant l'équation de Pell-Fermat $(F_d)$. 

\begin{corollaire} --- Soit $d\in\N$ un entier qui n'est pas un carré parfait. Alors, l'équation de Pell-Fermat $(F_d)$ admet une infinité de solutions. 
\end{corollaire}

\emph{Preuve} --- Notons $T$ la période minimale du développement de $\sqrt{d}$ et écrivons$\sqrt{d}=\left[a_0, \overline{a_1, a_2, ..., a_{T}}\right]$ ce développement. Comme $(v_n)$ est $T-$périodique, $v_T=v_0=1$. Si $T$ est pair alors $(-1)^Tv_T=1$ et donc $\mathscr{E}_1 \neq \varnothing$ ce qui assure que $(F_d)$ a une infinité de solutions. Si $T$ est impair, $(-1)^Tv_T=-1$ donc $\mathscr{E}_{-1} \neq \varnothing$ et, dans ce cas aussi, $(F_d)$ a une infinité de solutions. \hfill $\square$

\bigskip

\begin{rem} --- Dans le cas de $(F_d)$, le recherche de la solution fondamentale revient à déterminer le premier indice $N$ non nul tel que $v_N=1$ (avec la contrainte supplémentaire $N$ pair si $T$ est impair). La démonstration précédente montre que $N\leqslant T$. Or, en conservant les notations de la propriété \ref{prop01112157}, par définition, comme $N\geqslant 1$, $u_{N+1}=a_{N}v_{N}-u_{N}$ donc $a_{N}=u_{N+1}+u_{N}$. Or, par définition, $x_N=\frac{u_N+\sqrt{d}}{v_N}=u_N+\sqrt{d}$ donc $a_N=\entiere{x_N}=u_N+\entiere{\sqrt{d}}=u_N+a_0$. On en déduit que $u_{N+1}=a_0$. De plus, comme $v_N=1$, $v_N$ divise $u_N$ donc en utilisant les inégalités (\ref{eq20042001}) avec $j=0$, il vient $u_N=u_{N+1}$ et, ainsi, $a_N=2a_0$. Le point 5. de la proposition \ref{prop01112157} permet de conclure que $N=T$. Ainsi, le plus petit entier $N$ non nul tel que $v_N=1$ est $N=T$ et, par suite, $v_j=1$ si et seulement si $T$ divise $j$. On en déduit que la solution fondamentale de $(F_d)$ est donnée par la réduite $R_{T-1}$ si $T$ est pair et $R_{2T-1}$ si $T$ est impair.
\end{rem}

\subsection{\'Etude du cas général}

Nous atteignons (enfin!) notre but i.e. l'étude de l'équation diophantienne d'inconnue $(x,y)\in\N^2$

\[(E_2): ax^2-by^2=1\]

où $a$ et $b$ sont deux entiers naturels non nuls et premiers entre eux.

Remarquons que si $a=1$ ou $b=1$, on est dans le cas d'une équation de Pell-Fermat ($(F_{b}): x^2-by^2=1$ si $a=1$ et $(F_{a,-1})$: $y^2-ax^2=-1$ si $b=1$). On peut donc exclure ces cas dans la suite. 

On peut également exclure le cas trivial où $ab$ est un carré parfait. En effet, supposons qu'il existe un entier naturel non nul $k$ tel que $ab=k^2$ et que $(x,y)$ soit une solution de $(E_2)$. Comme $a$ et $b$ sont premiers entre eux, il existe deux entiers naturels non nuls $r$ et $t$ tels que $a=r^2$ et $b=t^2$. Alors, $(rx)^2-(ty)^2=1$ donc $(rx-ty)(rx+ty)=1$ ce qui impose que $rx-ty=rx+ty=1$. Il s'ensuit que $2rx=2$ i.e. $r=x=1$ et donc $y=0$ et $a=1$. Dans ce cas, l'équation $(E_2)$ est donc en fait une équation de Pell-Fermat qui n'admet qu'une solution triviale.

On peut donc se restreindre aux cas où ni $a$ ni $b$ n'est égal à 1 et où $ab$ n'est pas un carré parfait. 

Enfin, notons que la résolution dans $\N^2$ permet de déterminer, aux signes près, les solutions dans $\Z^2$ et que, comme nous ne faisons pas d'hypothèse sur l'ordre entre $a$ et $b$, la résolution de $(E_2)$ permet aussi celle de $ax^2-by^2=-1$.

\bigskip

\begin{theoreme} --- Soit $a$ et $b$ deux entiers supérieurs ou égaux à 2, premiers entre eux et tels que $d:=ab$ ne soit pas un carré parfait. On note $(u_n)$ et $(v_n)$ les suite associées à $\sqrt{d}$ par la propriété \ref{prop01112157}, $R_n=\dfrac{p_n}{q_n}$ la réduite d'indice $n$ dans le développement en fraction continuée de $\sqrt{d}$ et $T$ la plus petite période de ce développement. On note, de plus, $N=\frac{T}{2}$. 
\begin{enumerate}
\item Si $a<b$ alors l'équation $(E_2)$ admet une solution si et seulement si $T\equiv 0~[4]$, $v_N=a$ et $v_N$ divise $u_N$. De plus, dans ce cas, l'ensemble des solutions de $(E_2)$ est
\[\enstq{\left(\dfrac{p_{N-1+\ell T}}{a}, q_{N-1+\ell T}\right)}{\ell\in\N}.\]
\item Si $b<a$ alors l'équation $(E_2)$ admet une solution si et seulement si $T\equiv 2~[4]$, $v_N=b$ et $v_N$ divise $u_N$. De plus, dans ce cas, l'ensemble des solutions de $(E_2)$ est 
\[\enstq{\left(q_{N-1+\ell T}, \dfrac{p_{N-1+\ell T}}{b}\right)}{\ell\in\N}.\]
\item En particulier, $(E_2)$ n'a pas de solution si $T$ est impair.
\end{enumerate}
\label{theo04110947}
\end{theoreme}

\emph{Preuve}
\begin{enumerate}
\item Supposons $a<b$. Notons $(s,t)$ une solution de $ax^2-by^2=1$. En multipliant $as^2-bt^2=1$ par $a$, il vient $a^2s^2-abt^2=a$ donc $(as)^2-dt^2=a$ i.e. $(as,t)$ est solution de $(F_{d, a}):X^2-dY^2=a$ avec $1\leqslant a <\sqrt{d}$ car $a<b$. De plus, comme $s(as)-(bt)t=1$, le théorème de Bézout assure que $\text{PGCD}(as,t)=1$. Ainsi, d'après le théorème \ref{theo01112218}, il existe un entier $N\in\llbracket 1, T \rrbracket$ et un entier $\ell\in\N$ tel que $(-1)^Nv_N=a$ ou $(-1)^Nv_N=-a$ (ce deuxième cas ne pouvant se produire que si $T$ est impair) et tel que $(as, t)=(p_{N+\ell T-1},q_{N+\ell-1})$. Comme $v_N$ et $a$ sont positifs, on en déduit que $v_N=a$ et que $N$ est pair. Posons, pour simplifier les notations, $j=N+\ell T$ avec $j\geqslant 1$ car $N\geqslant 1$. Alors, comme $as=p_{j-1}$, $v_N=a$ divise $p_{j-1}$. Or, d'après (\ref{eq02111252}), $p_{j - 1}=q_{j - 1}u_{j} + q_{j - 2}v_{j}$ et, comme $(u_n)$ et $(v_n)$ sont $T-$périodiques, on a donc $p_{j-1}=q_{j-1}u_N+q_{j-2}v_N$. Ainsi, on en déduit que $v_N$ divise $q_{j- 1}u_{N}$. Mais, comme $p_{j - 1}$ et $q_{j- 1}$ sont premiers entre eux et $v_N$ divise $p_{j - 1}$, on peut affirmer que $v_N$ est premier avec $q_{j - 1}$ et le théorème de Gauss assure que $v_N$ divise $u_{N}$. Remarquons, de plus, que, $a$ étant différent de $1$, $v_N\neq 1$ donc $N \neq T$ puisque, par définition, $v_T=v_0=1$. Ainsi, $N\in\llbracket 1, T-1 \rrbracket$ et le point 5. de la propriété \ref{prop01112157} assure que $T=2N$ et donc, comme $N$ est pair, $T\equiv 0~[4]$. Enfin, $(s,t)$ est bien de la forme $\left(\dfrac{p_{N-1+\ell T}}{a}, q_{N-1+\ell T}\right)$.

Réciproquement, supposons que $T\equiv 0~[4]$, que $v_N=a$ et que $v_N$ divise $u_N$. Selon le point 5. de la propriété 6, $T=2N$ donc $N$ est pair et $a=(-1)^Nv_N$. Ainsi, d'après le théorème \ref{theo01112218}, pour tout entier $\ell\in\N$, $(p_{N+\ell T-1}, q_{N+\ell T-1})$ est une solution de $(F_{d,a})$ i.e. si on pose $j=N+\ell T$, $p_{j-1}^2-abq_{j-1}^2=a$ donc $a\left(\dfrac{p_{j-1}}{a}\right)^2-bq_{j-1}^2=1$. Il suffit alors pour conclure de montrer que $a$ divise $p_{j-1}$. Or, $j\geqslant 1$ car $N\geqslant 1$ donc en réutilisant le fait que, d'après (\ref{eq02111252}), 
\[p_{j-1}=q_{j-1}u_{j}+q_{j-2}v_{j}=q_{j-1}u_{N}+q_{j-2}v_N,\]
et l'hypothèse $v_N=a$ divise $u_N$, on peut affirmer que $a$ divise $p_{j-1}$. 

\item Supposons $b<a$. Notons $(s,t)$ une solution de $ax^2-by^2=1$. En multipliant $as^2-bt^2=1$ par $-b$, il vient $b^2t^2-abs^2=-b$ donc $(bt)^2-ds^2=-b$ i.e. $(bt,s)$ est solution de $(F_{d, -b}):X^2-dY^2=-b$ avec $1\leqslant \abs{-b} <\sqrt{d}$ car $b<a$ et $\text{PGCD}(bt,s)=1$. On est alors ramené à un cas similaire au précédent. En effet, il existe des entiers $N\in\llbracket 1, T\rrbracket$ et $\ell\in\N$ tels que $(-1)^Nv_N=-b$ ou $(-1)^Nv_ N=b$ (seulement possible si $T$ est impair) et $(bt,s)=(p_{N+\ell T-1}, q_{N+\ell T-1})$. On a alors $v_N=b$ car $v_N$ et $b$ sont positifs. En posant $j=N+\ell T$ avec $j\geqslant 1$, on peut raisonner comme dans le point 1. et en déduire que $T=2N$ En particulier, $T$ est pair donc on est dans le cas où $(-1)^Nv_N=-b$ ce qui impose que $N$ est impair. Il s'ensuit que $T\equiv 2~[4]$ et $(s,t)$ est bien de la forme $\left(q_{N-1+\ell T}, \dfrac{p_{N-1+\ell T}}{b},\right)$.

Réciproquement, supposons que $T\equiv 2~[4]$, que $v_N=b$ et que $v_N$ divise $u_N$. Alors, $T=2N$ donc $N$ est impair et $-b=(-1)^Nv_N$ donc, d'après le théorème \ref{theo01112218}, pour tout entier $\ell\in\N$, $(p_{N+\ell T-1}, q_{N+\ell T-1})$ est une solution de $(F_{d,-b})$ i.e. si on pose $j=N+\ell T$, $p_{j-1}^2-abq_{j-1}^2=-b$ donc $aq_{j-1}^2-b\left(\dfrac{p_{j-1}}{b}\right)^2=1$. On conclut en montrant, comme en 1., que $b$ divise $p_{j-1}$. \hfill $\square$ 
\end{enumerate}

\bigskip

Nous examinons à présent quelques exemples qui montrent que les trois conditions du théorème \ref{theo04110947} (reste de $T$ modulo 4, $v_N=a$ et $v_N$ divise $u_N$) sont minimales en ce sens que deux d'entre elles peuvent être réalisées sans que la troisième le soit. 

\bigskip

\begin{ex} \hspace{1cm}
\begin{enumerate}
\item Cas où $a<b$
\begin{enumerate}
\item Si $a=18$ et $b=23$ alors $d=414$ et la suite $(v_n)$ est $8-$périodique avec $u_4=v_4=18$ donc, comme $T=8\equiv 0~[4]$, les solutions de $18x^2-23y^2=1$ sont tous les couples de la forme $\left(\dfrac{p_{3+8k}}{18}, q_{3+8k}\right)$ avec $k\in\N$ où $\dfrac{p_n}{q_n}$ est la réduite d'indice $n$ dans le développement en fraction continuée de $\sqrt{414}$. Les 3 premières solutions sont $(26, 23)$, $(\np{1265394}, \np{1119433})$, $(\np{61586725954}, \np{54482804087})$.
\item Si $a=19$ et $b=25$ alors $d=475$ et la suite $(v_n)$ est $10-$périodique avec $v_5=u_5=19$. Ainsi, on a bien $v_5=a$ qui divise $u_5$ mais l'équation $19x^2-25y^2=1$ n'a pas de solution car $T=10\equiv 2~[4]$.
\item Si $a=18$ et $b=25$ alors $d=450$ et la suite $(v_n)$ est $8-$périodique avec $v_4=9$ et $u_4=18$. Ainsi, on a bien $T=8\equiv 0~[4]$ et $v_4$ qui divise $u_4$ mais l'équation $18x^2-25y^2=1$ n'a pas de solution car $v_4\neq a$.
\item Si $a=16$ et $b=19$ alors $d=304$ et la suite $(v_n)$ est $12-$périodique avec $v_6=16$ et $u_6=8$. Ainsi, on a bien $T=12\equiv 0~[4]$ et $v_6=a$ mais l'équation $16x^2-19y^2=1$ n'a pas de solution car $v_6$ ne divise pas $u_6$.
\end{enumerate} 
\item Cas où $b<a$
\begin{enumerate}
\item Si $a=25$ et $b=19$, d'après l'exemple 1. (b), $T=10\equiv 2~[4]$ et $u_5=v_5=b$ donc l'ensemble des solutions de l'équation $25x^2-19y^2=1$ est $\enstq{\left(q_{4+10\ell}, \dfrac{p_{4+10\ell}}{19}\right)}{\ell\in\N}$. Les 3 premières solutions sont $(34, 39)$, $(\np{3930298},\np{4508361})$ et $(\np{454334588170}, \np{521157514839})$
\item Si $a=23$ et $b=18$ alors, d'après l'exemple 1. (a), $T=8 \equiv 0~[4]$ donc $23x^2-18y^2=1$ n'a pas de solution.
\item 1. (c) et 1. (d) fournissent des contres-exemples similaires en échangeant $a$ et $b$.
\end{enumerate}
\end{enumerate}
\label{ex07062052}
\end{ex}

\bigskip

Nous terminons par examiner un cas favorable où la situation est plus simple.

\begin{corollaire} --- Soit $a$ et $b$ deux entiers premiers entre eux tels que $2 \leqslant a < b$  et tels que $d:=ab$ ne soit pas un carré parfait. On note $(u_n)$ et $(v_n)$ les suites associées à $\sqrt{d}$ par la propriété \ref{prop01112157} et $T$ la plus petite période du développement en fraction continuée de $\sqrt{d}$. On suppose que $T\equiv 0~[4]$ et que $a$ est impair. On pose $N=\frac{T}{2}$. Alors, l'équation $(E_2)$ admet des solutions si et seulement $v_N=a$. 
\label{coro07062020}
\end{corollaire}
 
\emph{Preuve} --- La condition $v_N=a$ est nécessaire d'après le théorème \ref{theo04110947}. Montrons que si $a$ est impair, elle est suffisante. On sait grâce à la propriété \ref{prop01112157} que le développement en fraction continuée de $\sqrt{d}$ est $\sqrt{d}=\left[a_0, \overline{a_1, a_2, ..., a_{T-1}, 2a_0}\right]$ où $(a_1, a_2, ..., a_{T-1})$ est un palindrome et celui de $x:=\entiere{\sqrt{d}}+\sqrt{d}$ est $x=\left[\overline{2a_0, a_1, ..., a_{T-1}}\right]$. Notons, pour tout $n\in\N$, $x_n$ le quotient complet d'indice $n$ de $x$. Alors, par définition, $x_N=\left[\overline{a_N, a_{N+1}, ..., a_{T-1}, 2a_0, a_1, ..., a_{N-2}, a_{N-1}}\right]$. De plus, d'après le corollaire \ref{coro18042015}, $x_N$ est un irrationnel quadratique réduit donc, d'après la propriété \ref{prop01111315}, $-\dfrac{1}{x_N^*}=\left[\overline{a_{N-1}, a_{N-2}, ..., a_1, 2a_0, a_{T-1} ..., a_{N+1}, a_{N}}\right]$. Mais alors, étant donné que $(a_1, a_2, ..., a_{T-1})$ est un palindrome et que $N=\frac{T}{2}$, 
\[-\dfrac{1}{x_N^*}=\left[\overline{a_{N+1}, a_{N+2}, ..., a_{T-1}, 2a_0, a_1 ..., a_{N-1}, a_{N}}\right]=x_{N+1}.\]
Par ailleurs, sachant que, d'après (\ref{eq20041953}), $v_N=\frac{d-u_N^2}{v_{N-1}}$,
\[-\dfrac{1}{x_N^*}=-\dfrac{1}{\frac{u_N-\sqrt{d}}{v_N}}=\dfrac{v_N}{\sqrt{d}-u_N}=\frac{v_N(\sqrt{d}+u_N)}{d-u_N^2}=\frac{u_N+\sqrt{d}}{v_{N-1}}.\]
Finalement, on aboutit à
\[\frac{u_N+\sqrt{d}}{v_{N-1}}=-\dfrac{1}{x_N^*}=x_{N+1}=\frac{u_{N+1}+\sqrt{d}}{v_{N+1}}\]
et donc, comme $(1, \sqrt{d})$ est une $\Q-$base de $\Q[\sqrt{d}]$, $v_{N+1}=v_{N-1}$ et $u_{N}=u_{N+1}$. 

Or, d'après (\ref{eq20041953}), $u_N=u_{N+1}=a_{N}v_N-u_N$ donc $v_N$ divise $2u_N$. Il s'ensuit que si $v_N=a$ est impair alors $v_N$ divise $u_N$ et donc, d'après le théorème \ref{theo04110947}, $(E_2)$ admet des solutions. \hfill $\square$

\bigskip

\begin{rem} \hspace{1cm}
\begin{enumerate}
\item Si $a$ est pair, l'équation peut ne pas avoir de solution même si $T\equiv 0~[4]$ et $v_N=a$ comme on l'a vu dans l'exemple \ref{ex07062052} avec $a=16$ et $b=19$.
\item On peut donner un résultat similaire si $b<a$ en supposant que $T\equiv 2~[4]$ et que $b$ est impair.
\end{enumerate}

\end{rem} 

\bigskip

Signalons, enfin, que K. Matthews propose à l'adresse
\begin{center}
\texttt{http\string://www.numbertheory.org/php/pellab.html}
\end{center}
un programme permettant de déterminer si l'équation $(E_2)$ admet ou non des solutions et qui détermine la première solution. L'algorithme est basé sur un théorème donné dans la note très succincte \cite{Mat12}. La démarche est sensiblement la même que la nôtre mais en considérant les quotients complets de $\sqrt{\frac{b}{a}}$ ce qui fournit une condition d'existence de solutions un peu plus simple que celle donnée dans le théorème 2 ci-dessus. Matthews ne donne cependant pas de description effective de l'ensemble des solutions.

\section{Le cas $\boldsymbol{n\geqslant 3}$}
\label{secThue}
\begin{theoreme} --- Soit $a$ et $b$ deux entiers relatifs non nuls et $n$ un entier naturel non nul. Si $n\geqslant 3$ alors l'équation $(E_n): ax^n-by^n=1$ d'inconnue $(x,y)\in\Z^2$ n'a qu'un nombre fini de solutions.
\label{th03051542}
\end{theoreme}

\noindent On déduit ce résultat d'un théorème démontré par A. Thue en 1909 (\cite{Thu09}, \cite[p. 120]{Du07}, \cite[p. 231]{Hi08}).

\bigskip

On rappelle qu'un nombre complexe $\alpha$ est dit algébrique s'il est racine d'un polynôme non nul de $\Q[X]$. Dans ce cas, l'ensemble $\mathcal{I}_{\alpha}:=\enstq{P\in\Q[X]}{P(\alpha)=0}$ est un idéal non réduit à $\{0\}$ de $\Q[X]$. Il est donc monogène et l'unique polynôme unitaire $\mu_{\alpha}$ qui engendre $\mathcal{I}_{\alpha}$ est appelé le polynôme minimal de $\alpha$. De plus, si on note $d$ le degré de $\mu_{\alpha}$, on dit que $\alpha$ est un nombre algébrique de degré $d$.

\begin{theoreme} (Thue, 1909) --- Si $\alpha$ est un nombre algébrique de degré $d\geqslant 2$ alors, pour tout $\varepsilon>0$, il existe une constante $C>0$ (dépendant de $\alpha$ et $\varepsilon$) telle que, pour tout nombre rationnel $\frac{p}{q}$ (avec $p\in\Z$ et $q\in\N^*$),
\[\abs{\alpha - \frac{p}{q}} \geqslant \frac{C}{q^{\frac{d}{2}+1+\varepsilon}}.\] 
\label{thThue}
\end{theoreme}

\emph{Démonstration du théorème \ref{th03051542}}

Soit $n$ un entier au moins égal à $3$ et $(s,t)\in\Z^2$ une solution de $(E_n)$.

Commençons par remarquer que si $s=0$ alors $t$ est solution de $y^n=-\frac{1}{b}$, ce qui n'est possible que si $b=\pm1$ et alors $t\in\{-1, 1\}$. De même, si $t=0$ alors $s$ est solution de $x^n=\frac{1}{a}$ donc $s\in\{-1, 1\}$. 

Il suffit donc de démontrer que $(E_n)$ admet un nombre fini de solutions dont les deux termes sont non nuls. Dans la suite, on suppose $st\neq 0$.

 Puisque $as^n=bt^n+1$, $as^n$ et $bt^n$ sont des entiers non nuls consécutifs  donc $\left|{\abs{as^n}-\abs{bt^n}}\right|=1$.

Notons $\alpha=\sqrt[n]{\frac{\abs{a}}{\abs{b}}}$.

\medskip

\emph{1\up{er} cas}. --- Supposons que $\alpha$ soit rationnel. Alors, il existe deux entiers naturels non nuls et premiers entre eux $u$ et $v$ tels que $\alpha=\frac{u}{v}$. Dès lors, $\frac{\abs{a}}{\abs{b}}=\frac{u^n}{v^n}$ et, comme ces deux fractions sont irréductibles, $\abs{a}=u^n$ et $\abs{b}=v^n$. Ainsi, $\abs{\abs{u^ns^n}-\abs{v^nt^n}}=1$ i.e. $\abs{\abs{us}^n-\abs{vt}^n}=1$. Or, comme $\abs{us}\geqslant 1$ et $\abs{vt}\geqslant 1$,
\[\abs{\abs{us}^n-\abs{vt}^n}=\abs{\abs{us}-\abs{vt}}\sum_{k=0}^{n-1}\abs{us}^k\abs{vt}^{n-1-k} \geqslant 1\times n \geqslant 3.\]
On aboutit à une contradiction donc, si $\alpha$ est rationnel, $(E_n)$ n'a pas de solution dont les deux termes sont non nuls.

\medskip

\emph{2\up{e} cas}. --- Supposons que $\alpha$ ne soit pas un nombre rationnel. Comme $\alpha$ est racine du polynôme $P:=\abs{b}X^n-\abs{a}$, $\alpha$ est un nombre algébrique de degré $d\in\llbracket 2, n\rrbracket$. En appliquant le théorème de Thue avec $\varepsilon=\frac{1}{4}$, on en déduit qu'il existe un réel $C>0$ indépendant de $s$ tel que  
\begin{equation}
\abs{\alpha - \frac{\abs{t}}{\abs{s}}}\geqslant \frac{C}{\abs{s}^{\frac{2d+5}{4}}}.
\label{eq04051836}
\end{equation}
Or, en divisant $\abs{\abs{as^n}-\abs{bt^n}}=1$ par $\abs{bs^n}$, il vient $\abs{\frac{\abs{a}}{\abs{b}}-\frac{\abs{t}^n}{\abs{s}^n}}=\frac{1}{\abs{bs^n}}$ i.e. $\abs{\alpha^n-\left(\frac{\abs{t}}{\abs{s}}\right)^n}=\frac{1}{\abs{bs^n}}$. On en déduit que
\begin{equation}
\frac{1}{\abs{bs^n}}=\abs{\alpha-\frac{\abs{t}}{\abs{s}}}\sum_{k=0}^{n-1}\alpha^k\left(\frac{\abs{t}}{\abs{s}}\right)^{n-k-1} \geqslant \abs{\alpha-\frac{\abs{t}}{\abs{s}}}\alpha^{n-1}.
\label{eq04051838}
\end{equation}
car tous les termes de la somme sont positifs dont on peut la minorer par son dernier terme. 

On déduit alors de (\ref{eq04051836}) et (\ref{eq04051838}) que $\dfrac{1}{\abs{bs^n}}\geqslant \dfrac{\alpha^{n-1}C}{\abs{s}^{\frac{2d+5}{4}}}$ et donc $\abs{s}^{\frac{4n-2d-5}{4}} \leqslant \dfrac{1}{\alpha^{n-1}\abs{b}C}$. 

Or, puisque $d\leqslant n$ et $n\geqslant 3$, $\dfrac{4n-2d-5}{4}\geqslant \dfrac{2n-5}{4} \geqslant \dfrac{1}{4}$. 

Ainsi, comme $\abs{s}\geqslant 1$, $\abs{s}^{\frac{1}{4}} \leqslant \dfrac{1}{\alpha^{n-1}\abs{b}C}$ et donc $\abs{s} \leqslant \dfrac{1}{(\alpha^{n-1}\abs{b}C)^4}$.

Il s'ensuit que $s$ ne peut prendre qu'un nombre fini de valeurs et, comme pour une valeur de $s$ donnée, il y, au plus, 2 valeurs possibles pour $t$, on en déduit que $(E_n)$ a nombre fini de solutions. \hfill $\square$

\bigskip

\bigskip

\noindent \textbf{Remerciements} --- L'auteur remercie chaleureusement Daniel Duverney, Alain Faisant et Marc Hindry pour leurs réponses, conseils et encouragements dans la rédaction de cet article ainsi que B. Randé et G. Alarcon du comité de lecture de la RMS pour leurs relectures attentives et pour les corrections proposées.

\nocite{Hi08}
\nocite{HW60}
\nocite{Mo69}
\nocite{NZM}
\nocite{Hu82}
\nocite{Mo01}
\nocite{Mol04}
\nocite{Mat12}
\nocite{RS94}
\nocite{Fai91}

\bibliographystyle{smfalpha}
\bibliography{biblio} 

\providecommand{\bysame}{\leavevmode ---\ }
\providecommand{\og}{``}
\providecommand{\fg}{''}
\providecommand{\smfandname}{et}
\providecommand{\smfedsname}{\'eds.}
\providecommand{\smfedname}{\'ed.}
\providecommand{\smfmastersthesisname}{M\'emoire}
\providecommand{\smfphdthesisname}{Th\`ese}
\begin{thebibliography}{{Whi}12}

\bibitem[Ben01]{Be01}
{\scshape M.~A. Bennett} -- {\og Rational approximation to algebraic numbers of
  small height: the {D}iophantine equation $\abs{ax^n-by^n}=1$\fg}, \emph{J.
  Reine Angew. Math.} \textbf{525} (2001), p.~1--49.

\bibitem[Duv07]{Du07}
{\scshape D.~Duverney} -- \emph{Théorie des nombres}, 2nde \smfedname, Dunod,
  2007.

\bibitem[{Fai}91]{Fai91}
{\scshape A.~{Faisant}} -- \emph{{L'\'equation diophantienne du second
  degr\'e}}, Hermann, {P}aris, 1991.

\bibitem[Hin08]{Hi08}
{\scshape M.~Hindry} -- \emph{Arithmétique}, Tableau noir, Calvage \& Mounet,
  2008.

\bibitem[Hua82]{Hu82}
{\scshape L.~K. Hua} -- \emph{Introduction to number theory}, Spinger-Verlag,
  1982.

\bibitem[HW60]{HW60}
{\scshape G.~H. Hardy {\normalfont \smfandname} E.~M. Wright} -- \emph{An
  {I}ntroduction to the {T}heory of {N}umbers}, 6$^\text{e}$ \smfedname, Oxford
  University Press, 1960.

\bibitem[LJ02]{Le02}
{\scshape H.~Lenstra~Jr} -- {\og Solving the {P}ell equation\fg}, \emph{Notices
  Am. Math. Soc.} \textbf{49} (2002), p.~182--192.

\bibitem[{Mat}12]{Mat12}
{\scshape K.~{Matthews}} -- {\og {A {M}idpoint {C}riterion for the
  {D}iophantine {E}quation $ax^2 - by^2 = \pm1$}\fg}, disponible à l'adresse
  \url{http://www.numbertheory.org/pdfs/keith_midpoint.pdf}, 2012.

\bibitem[Mol01]{Mo01}
{\scshape R.~A. Mollin} -- {\og Quadratic diophantine equations determined by
  continued fractions\fg}, \emph{{JP} {J}our. {A}lgebra, {N}umber {T}heory \&
  {A}ppl.} \textbf{1} (2001), p.~57--75.

\bibitem[{Mol}04]{Mol04}
{\scshape R.~A. {Mollin}} -- {\og {A continued fraction approach to the
  {D}iophantine equation $ax^2- by^2= \pm1$}\fg}, \emph{{JP J. Algebra Number
  Theory Appl.}} \textbf{4} (2004), p.~159--207.

\bibitem[Mor69]{Mo69}
{\scshape L.~J. Mordell} -- \emph{Diophantine {E}quations}, Academic {P}ress,
  1969.

\bibitem[NZM91]{NZM}
{\scshape I.~Niven, H.~S. Zuckerman {\normalfont \smfandname} H.~L. Montgomery}
  -- \emph{An {I}ntroduction to the {T}heory of {N}umbers}, 5$^\text{e}$
  \smfedname, John Wiley \& Sons, Inc., 1991.

\bibitem[RS94]{RS94}
{\scshape A.~M. {Rockett} {\normalfont \smfandname} P.~{Sz\"usz}} --
  \emph{{Continued fractions}}, World Scientific, 1994.

\bibitem[{Thu}09]{Thu09}
{\scshape A.~{Thue}} -- {\og {\"Uber {A}nn\"aherungswerte algebraischer
  {Z}ahlen}\fg}, \emph{{J. {R}eine {A}ngew. {M}ath}} \textbf{135} (1909),
  p.~284--305.

\bibitem[{Whi}12]{Whi12}
{\scshape E.~E. {Whitford}} -- \emph{{The {P}ell equation}}, New York, 1912.

\end{thebibliography}

\end{document}